\documentclass[11pt]{article}

\usepackage{latexsym}
\usepackage[all]{xy}
\usepackage{amsmath}
\usepackage{amssymb}
\usepackage{amsfonts}
\usepackage{accents}

\newcommand  {\dgcat}     {dg-\mathcal{C}at}

\newcommand  {\dAff}     {\mathbf{dAff}}

\newcommand  {\dSt}   {\mathbf{dSt}}
\newcommand  {\hh}   {\mathbf{H}}
\newcommand {\s}{\infty}
\newcommand  {\Sch}   {\mathbf{Sch}}
\newcommand  {\dSch}   {\mathbf{dSch}}

\newtheorem{thm}{Theorem}[section]
\newtheorem{prop}[thm]{Proposition}
\newtheorem{lem}[thm]{Lemma}
\newtheorem{sublem}[thm]{Sublemma}
\newtheorem{df}[thm]{Definition}

\newtheorem{conj}[thm]{Conjecture}

\begin{document}

\title{\textbf{Proper local complete intersection morphisms preserve perfect complexes}}  
\author{\bigskip\\
Bertrand To\"en\thanks{This work is partially supported by the ANR grant ANR-09-BLAN-0151 (HODAG).}
 \\
\small{I3M, Universit\'e de Montpellier 2}\\
\small{Case Courrier 051}\\
\small{Place Eug\`ene Bataillon}\\
\small{34095 Montpellier Cedex, France}\\
\small{e-mail: btoen@math.univ-montp2.fr}}

\date{October 2012}

\maketitle

\begin{abstract}
Let $f : X \longrightarrow Y$ be a proper and local complete intersection morphism of schemes. We prove
that $\mathbb{R}f_{*}$ preserves perfect complexes, without any projectivity or noetherian assumptions. 
This provides a different proof of a theorem by Neeman and Lipman (see \cite{neli}) 
based on techniques from derived algebraic geometry to proceed a reduction to the noetherian case. 
\end{abstract}

\section*{Introduction}

In \cite[Exp. XIV]{sga6} Grothendieck and its collaborators present a huge list of open problems around
the most general version of intersection theory on schemes and of the Grothendieck-Riemann-Roch formula. 
An important part of these problems concern the question of proving the GRR formula for proper
morphisms of schemes without assuming the existence of a global factorization and/or without 
any noetherian assumption. One of the most basic questions is the existence of
push-forward operations in algebraic K-theory of schemes in the most general setting. Directly related to this
specific question is the following conjecture.

\begin{conj}(See \cite[ExpIII-2.1]{sga6})\label{conj}
Let $f : X \longrightarrow Y$ be a proper and pseudo-coherent morphism of schemes. Then, the derived direct
image
$$\mathbb{R}f_{*} : D_{qcoh}(X) \longrightarrow D_{qcoh}(Y)$$
preserves pseudo-coherent complexes. 
\end{conj}

This conjecture has been proven in \cite{ki}, and more recently A. Neeman and J. Lipman deduced
from it the following theorem.

\begin{thm}(\cite[Ex. 2.2 (a)]{neli})
Let $f : X \longrightarrow Y$ be a proper and perfect morphism of schemes. 
Then, the derived direct image
$$\mathbb{R}f_{*} : D_{qcoh}(X) \longrightarrow D_{qcoh}(Y)$$
preserves perfect complexes. 
\end{thm}

The purpose of this work is to propose a new proof of the above theorem by assuming further that
 $f$ is a local complete intersection morphism. 

\begin{thm}\label{ti}
Let $f : X \longrightarrow Y$ be a proper and local complete intersection morphism of schemes. 
Then, the derived direct image
$$\mathbb{R}f_{*} : D_{qcoh}(X) \longrightarrow D_{qcoh}(Y)$$
preserves perfect complexes. 
\end{thm}

The proof we propose in this work is based on a reduction to the noetherian case (no surprise) but in the somehow
unexpected context of derived algebraic geometry. The general strategy proceeds as follows. 
The statement being local on $Y$ we can reduce it to the case where $Y$ is an affine scheme $Spec\, A$. We
then write $A$ as a filtered colimit of noetherian rings $A_{i}$ and we try to descend the whole situation 
over one of the ring $A_{i}$. The unexpected fact is that even thought the scheme $X$ descend as a scheme $X_{i}$
proper and of finite type over some $A_{i}$, it does not seem to be true that $X_{i}$ can be chosen
to be itself a local complete intersection. Our main observation is that $X_{i}$ can be
chosen to be a proper and local complete intersection \emph{derived scheme}. This fact can then be used
to prove the theorem \ref{ti} by using some standard facts about derived categories of derived schemes, particularly
the base change property (see \ref{p1}), which is a statement specifically true in the derived setting 
and wrong in the underived setting without extra flatness assumptions. This method of proof also 
shows that the theorem \ref{ti} remains true for $X$ and $Y$ being themselves derived schemes, even thought we
do not make this statement explicit in this work. \\

\bigskip

\textbf{Acknowledgements:} I am grateful to G. Vezzosi, M. Vaqui\'e for their comments on a preliminary version of 
this work. I would like also to thank B. Duma, I have learned about the finiteness conjecture 
\cite[ExpIII-2.1]{sga6} while reading his thesis \cite{du}. I am very grateful to 
A. Neeman for brought to me the papers \cite{ki} and \cite{neli}, which I was not aware
before writting the very first version of this work.

\section{Derived schemes}

We have collected in this preliminary section some facts concerning derived schemes that we will
use in our proof of the theorem \ref{ti}. They belong to the general properties of derived schemes and
derived stacks and are certainly well known to experts. Some of these statements are probably established in 
the topological setting in \cite[\S 2]{lu}. We note however some differences between the theory of derived schemes
we will be using and the theory of spectral schemes used in \cite{lu} (notably cotangent complexes are not 
quite the same). We have therefore included proofs for the three
main statements we will be using: base change, continuity and noetherian approximation. \\

We recall from \cite{hagII,seat,tova} the existence of an $\s$-category $\dSt$, of derived stacks, and its full 
subcategory 
$\dSch \subset \dSt$ of derived schemes. The $\s$-category $\dSch$ contains as a full sub-category
$\Sch$, the category of schemes. The inclusion functor $\Sch \hookrightarrow \dSch$ possesses
a right adjoint $h^{0} : \dSch \longrightarrow \Sch$, the truncation functor. We note that 
even thought both $\s$-categories $\Sch$ and $\dSch$ admit finite limits, 
$\Sch$ is not stable by finite limits in $\dSch$.

Affine derived schemes are of the form $Spec\, A$ for some commutative
simplicial ring $A \in sComm$. The $Spec$ $\s$-functor defines an equivalence of $\s$-categories
$$Spec : \mathbf{sComm}^{op} \simeq\dAff,$$
where $\dAff \subset \dSch$ is the full sub $\s$-category of derived affine schemes, and
$\mathbf{sComm}^{op}$ is the $\s$-categories of simplicial commutative rings. Restricted
to affine objects the inclusion $\s$-functor from schemes to derived schemes, and its right adjoint
$h^{0}$, is equivalent to the adjunction of $\s$-categories
$$\pi_{0} : \mathbf{sComm} \longleftrightarrow \mathbf{Comm} : i,$$
where $i$ sends a ring to the constant simplicial ring and $\pi_{0}$ is the
connected component $\s$-functor. 

For all derived scheme $X$, the scheme $h^{0}(X)$, considered as an object in $\dSch$, has an adjuntion
morphism $j : h^{0}(X) \longrightarrow X$ which is a closed immersion. When $X=Spec\, A$ is affine, 
this morphism is equivalent to 
$Spec\, \pi_{0}(A) \longrightarrow Spec\, A$ corresponding to the natural projection $A \longrightarrow \pi_{0}(A)$. 
The higher homotopy groups $\pi_{i}(A)$ are $\pi_{0}(A)$-modules, and thus define quasi-coherent 
sheaves on $Spec\, \pi_{0}(A)$. These sheaves will be denoted by $h^{-i}(X):=\pi_{i}(A)$. This construction
glue in the non-affine case: for any derived scheme $X$ the scheme $h^{0}(X)$ carries natural
quasi-coherent sheaves $h^{i}(X):=\pi_{i}(\mathcal{O}_{X})$.

Any derived scheme $X$ possesses a cotangent complex $\mathbb{L}_{X} \in L_{qcoh}(X)$ 
(see \cite[\S 1.4]{hagII}, \cite[\S 4.2]{seat}), 
which is a quasi-coherent 
complex on $X$ (see our next subsection for quasi-coherent complexes and the definition of $L_{qcoh}(X)$).
When $X=Spec\, A$ is affine, $\mathbb{L}_{X}$ is the cotangent complex of $A$ introduced by 
Quillen (it is a simplicial $A$-module than can be turned into a dg-module over the normalization 
of $A$ to consider it as an object in $L_{qcoh}(X)$). For a morphism of derived schemes
$f : X \longrightarrow Y$ we define the relative cotangent complex $\mathbb{L}_{f}$, or 
$\mathbb{L}_{X/Y}$, as the cone, in the dg-category $L_{qcoh}(X)$, of the natural morphism
$$f^{*}(\mathbb{L}_{Y}) \longrightarrow \mathbb{L}_{X}.$$

\begin{df}\label{d0}
\begin{enumerate}

\item  Let $X$ be a derived scheme. We say that $X$ is \emph{quasi-compact} (resp. \emph{quasi-separated}, resp.
\emph{separated}) if the scheme $h^{0}(X)$ is so.

\item Let $f : X \longrightarrow $ be a morphism of derived schemes. We say that 
$f$ is \emph{proper} if $h^{0}(f) : h^{0}(X) \longrightarrow h^{0}(Y)$ is a proper morphism
of schemes. 

\item Let $f : X \longrightarrow Y$ be a morphism of derived schemes. We say that 
$f$ is \emph{locally of finite presentation} if $h^{0}(f) : h^{0}(X) \longrightarrow h^{0}(Y)$ 
is a morphism locally of finite presentation of schemes and if the
relative cotangent complex $\mathbb{L}_{f}$ is perfect on $X$.
\end{enumerate}
\end{df}

We note that the definition above of morphism locally of finite presentation possesses several
equivalent versions (see \cite[Prop. 2.2.2.4]{hagII}). In particular, a morphism of derived affine schemes
$f : Spec\, A \longrightarrow Spec\, B$ is locally of finite presentation if and only if
$B$ is equivalent to a retract of a finite cell commutative simplicial $B$-algebra (see 
\cite[Prop. 1.2.3.5]{hagII}, or \cite[Prop. 2.2]{tova}). 

We will also be using the following notion. 

\begin{df}\label{d01}
Let $X=Colim\, X_{\alpha}$ be a colimit in the $\s$-category $\dSch$ of derived schemes. We will
say that the colimit is \emph{strong} if it is also a colimit in the bigger $\s$-category
$\dSt$ of derived stacks. 
\end{df}

\subsection{Review of derived categories of derived schemes}

According to \cite[\S 2]{to2}, 
for any derived stack $X \in \dSt$, we have a natural ($\mathbb{Z}$-linear) dg-category 
$L_{qcoh}(X)$ of quasi-coherent complexes over $X$. This can be made into an $\s$-functor
$$L_{qcoh} : \dSt^{op} \longrightarrow \dgcat,$$
where $\dgcat$ is the $\s$-category of (locally presentable, see \cite{to2}) dg-categories
For a morphism of derived stacks
$f : X \longrightarrow Y$, we have an adjunction of dg-categories
$$f^{*} : L_{qcoh}(Y) \longleftrightarrow L_{qcoh}(X) : f_{*}.$$

When $X=Spec\, A$ is an affine derived scheme, the dg-category $L_{qcoh}(X)$ can be
explicitely described, up to a natural equivalence, as follows. The commutative
simplicial ring $A$ possesses a normalisation $N(A)$, which is a (commutative) dg-ring. 
The dg-category $L(A)$ of cofibrant $N(A)$-dg-modules is naturally equivalent to $L_{qcoh}(X)$. 
This description is moreover functorial in the following way. A morphism of derived
affine schemes $f : X=Spec\, A \longrightarrow Y=Spec\, B$ corresponds to a
morphism of commutative simplicial rings $B \longrightarrow A$, and thus to a morphism
of dg-rings $N(B) \longrightarrow N(A)$. The adjunction
$$f^{*} : L_{qcoh}(Y) \longleftrightarrow L_{qcoh}(X) : f_{*}$$
is then equivalent to the following adjunction
$$ -\otimes_{A}B : L(A) \longleftrightarrow L(B) : f,$$
where we have written $\otimes_{A}$ symbolically for $\otimes_{N(A)}N(B)$, 
and where $f$ is the forgetful functor. 

In the general case, we write a general derived stack $X$ as a colimit in $\dSt$,
$X=Colim\, Spec\, A_{\alpha}$, and we have
$$L_{qcoh}(X) \simeq Lim\, L_{qcoh}(Spec\, A_{\alpha}) \simeq Lim\, L(A_{\alpha}),$$
where the limit here is taken in the $\s$-category of dg-categories. When $X$ is
a derived scheme we can take a rather simple colimit description by taking all
the $Spec\, A_{\alpha}$ to belong to a basis of opens for the Zariski topology on $X$.

Finally, the dg-category $L_{qcoh}(X)$ of any derived scheme $X$ possesses a natural 
non-degenerate t-structure (by which we mean that the associated triangulated
category $[L_{qcoh}(X)]$ has such a t-structure), whose heart is canonically equivalent to the
category of quasi-coherent sheaves on the scheme $h^{0}(X)$. Locally on $X$, over an 
affine open $U=Spec\, A \subset X$ this t-structure can be described as follows. The dg-category 
$L_{qcoh}(U)$ is identified with $L(A)$ the dg-category of cofibrant $N(A)$-dg-modules. An object
$E$ in $L(A)$ is declared to belong to $L(A)^{\leq 0}$ if it is such that 
$H^{i}(E)=0$ for all $i\geq 0$. It is easy to see that this defines an aisle of a non-degenerate
t-structure on $D(A)=[L(A)]$, the derived category of $N(A)$-dg-modules. The heart is
the full sub-category of $D(A)$ consisting of dg-module $E$ with $H^{i}(E)=0$ except for $i=0$.
This sub-category is equivalent, via the functor $E \mapsto H^{0}(E)$, to the category 
of $\pi_{0}(A)$-modules. 

\begin{df}\label{d1}
Let $X$ be a derived stack. An object $E \in L_{qcoh}(X)$ is 
\emph{perfect} if for all affine derived scheme $Z=Spec\, A$ and all
morphism $f : Z \longrightarrow X$, the pull-back $f^{*}(E)$ is a
compact object in $L_{qcoh}(Z)$. 

The full sub-dg-category of $L_{qcoh}(X)$ consisting of perfect complexes
will be denoted by $L_{parf}(X)$.
\end{df}

Note that for $Z=Spec\, A$ affine, the compact objects in $L_{qcoh}(Z)$ can also be described 
in several different ways: as strongly dualizable objects or as retracts of finite
cell $A$-dg-modules (see discussion after \cite[Def. 2.3]{tova}). \\

Suppose now that we have a cartesian square of derived stacks
$$\xymatrix{
X' \ar[r]^-{g} \ar[d]_-{q} & Y' \ar[d]^-{p} \\
X \ar[r]_{f} & Y,}$$
by adjunction there is, for any object $E \in L_{qcoh}(Y')$, a natural morphism
in the dg-category $L_{qcoh}(X)$ 
$$\rho_{E} : f^{*}p_{*}(E) \longrightarrow q_{*}g^{*}(E).$$

\begin{prop}\label{p1}
If 
$$\xymatrix{
X' \ar[r]^-{g} \ar[d]_-{q} & Y' \ar[d]^-{p} \\
X \ar[r]_{f} & Y,}$$
is a cartesian square of quasi-compact and quasi-separated derived schemes, then 
for all $E \in L_{qcoh}(Y')$ the morphism $\rho_{E} : f^{*}p_{*}(E) \longrightarrow q_{*}g^{*}(E)$
above is an equivalence.
\end{prop}

\textit{Proof:} Localising on the Zariski topology on $X$ we can assume that $X$ and $Y$ are affine
derived schemes. Let us write the cartesian square as follows
$$\xymatrix{
Z_{B}:=Z\times_{Spec\, A}Spec\, B \ar[r]^-{g} \ar[d]_-{q} & Y'=Z \ar[d]^-{p} \\
X=Spec\, B \ar[r]_{f} & Y=Spec\, A.}$$

For $E \in L_{qcoh}(Z)$, the object $f^{*}p_{*}(E)$ is $\hh(Z,E)\otimes_{A}^{\mathbb{L}}B$, where
$\hh(Z,E)$ is the cohomology $A$-dg-module of cohomology of $Z$ with coefficients in $E$. The object
$q_{*}g^{*}(E)$ is $\hh(Z_{B},g^{*}(E))$, the cohomology $B$-dg-module of $Z_{B}$ with coefficients
in $g^{*}(E)$. The morphism $\rho_{E}$ is then the natural morphism $g^{*} : \hh(Z,E) \longrightarrow
\hh(Z_{B},g^{*}(E))$, extended to $\hh(Z,E)\otimes_{A}^{\mathbb{L}}B$ by linearity. 

As $Z$ is quasi-compact and quasi-separated it belongs to the smallest full sub $\s$-category 
of $\dSch$ containing affines and which is stable by finite strong colimits in $\dSt$ (see definition
\ref{d01}). Therefore, to prove 
our proposition it is enough to prove the following two individual statements:

\begin{enumerate}

\item If $Z$ is affine, then for all $E \in L_{qcoh}(Z)$, the morphism 
$\rho_{E}$ is an equivalence.

\item The full sub $\s$-category of objects $Z \in \dSch/Spec\, A$ for which $\rho_{E}$
is an equivalence for all $E \in L_{qcoh}(Z)$ is stable by finite strong colimits. 

\end{enumerate}

The property $(1)$ follows directly from the shape of fiber products of derived
affine schemes (see \cite[Prop. 1.1.0.8]{hagII}). For the property $(2)$, let 
$Z \simeq Colim\,  Z_{i}$ be a finite strong colimit in $\dSch$ 
for which we know that the proprosition holds for all the derived schemes $Z_{i}$. 
We let $Z_{i,B}=Z_{i}\times_{Spec\, A}Spec\, B$, and we notice that as colimits are
universal in $\dSt$ (because it is an $\s$-topos, see \cite{hagI,hagII,seat}), we have an induced strong colimit
$Z_{B} \simeq Colim \, Z_{i,B}$. 
As the colimit is strong we have moreover
$$\hh(Z,E)\simeq Lim\,  \hh(Z_{i},E_{|Z_{i}}),$$
where $Lim$ stands for the limit in the dg-category $L_{qcoh}(Spec\, A)$. As this limit is finite, we
have
$$\hh(Z,E)\otimes_{A}^{\mathbb{L}}B \simeq  Lim\, (\hh(Z_{i},E_{|Z_{i}})\otimes_{A}^{\mathbb{L}}B) \simeq
 Lim\, \hh(Z_{i,B},g^{*}(E)_{|Z_{i}}) \simeq\hh(Z_{B},g^{*}(B)).$$
\hfill $\Box$ \\

The second formal property we will need is continuity, relating the dg-category of 
quasi-coherent complexes of certain limits of derived schemes to the colimit (inside the $\s$-category 
of dg-categories) of the dg-categories of quasi-coherent complexes on each individual derived schemes. 

For this we let $A$ be a commutative ring wich is written as a filtered colimit
$$A = Colim\, A_{i}.$$
We will suppose that the indexing category $I$ has an initial obect $0 \in I$. 
Let $X_{0} \longrightarrow Spec\, A_{0}$ be a derived scheme, and let set 
$X_{i}:=X_{0}\times_{Spec\, A_{0}}Spec\, A_{i}$ its base change to $A_{i}$. We will also 
denote by $X:=X\times_{Spec\, A_{0}}Spec\, A$ the base change to $A$ itself. In this situation 
we have a natural morphism of dg-categories
$$Colim\, L_{qcoh}(X_{i}) \longrightarrow L_{qcoh}(X).$$

\begin{prop}\label{p2}
With the same notations as above, if the derived scheme $X_{0}$ is quasi-compact and quasi-separated
then the morphism
$$Colim\,  L_{parf}(X_{i}) \longrightarrow L_{parf}(X)$$
is an equivalence of dg-categories.
\end{prop}

\textit{Proof:} It follows the same lines as the proof of the proposition \ref{p1}. 
We reduce the proposition to the following two individual statements.

\begin{enumerate}

\item Proposition \ref{p2} holds for $X_{0}$ affine.

\item The full sub $\s$-category of objects $X_{0} \in \dSch/Spec\, A_{0}$ for which 
the proposition holds is stable by strong finite colimits.

The property $(1)$ simply states that for a filtered colimit of commutative
simplicial rings $B=Colim\, B_{i}$, we have
$$Colim\, L_{parf}(Spec\, B_{i}) \simeq L_{parf}(Spec\, B),$$
which is a particular case of \cite[Lem. 2.10]{tova}. The second property is proven as follows. 
Let $X_{0} = Colim\, X_{0,\alpha}$ be a finite strong colimit in $\dSch$ such that
the proposition \ref{p2} holds for all $X_{0,\alpha}$. We let 
$X_{i,\alpha}:=X_{i}\times_{X_{0}}X_{0,\alpha}$ and $X_{\alpha}=X\times_{X_{0}}X_{0,\alpha}$.
As the colimit is filtered, it commutes with finite limits, and thus we have
$$Colim_{i} \, L_{parf}(X_{i}) \simeq Colim_{i}\,  Lim_{\alpha} \, L_{parf}(X_{i,\alpha}) \simeq
Lim_{\alpha} \, Colim_{i} \, L_{parf}(X_{i,\alpha}) \simeq Lim_{\alpha} \, L_{parf}(X_{\alpha}).$$
\hfil $\Box$ \\

\end{enumerate}

\subsection{Noetherian approximation for derived schemes}

We let $A$ be a commutative ring which is written as a filtered colimit $A=Colim \, A_{i}$. 
The $\s$-category of derived schemes over $Spec\, A_{i}$ (resp. over $Spec\, A$) will be denoted by 
$\dSch_{A_{i}}$ (resp. $\dSch_{A}$). 
We will study the $\s$-functor
$$Colim\,  \dSch_{A_{i}} \longrightarrow \dSch_{A}.$$

For this we denote by $\dSch^{\leq n}_{A}$ the full sub $\s$-category of
$\dSch_{A}$ consisting of derived schemes $f : X \longrightarrow Spec\, A$, 
for which $f$ is of locally of finite presentation, quasi-compact, quasi-separated, and such that
the cotangent complex $\mathbb{L}_{f}$ is of amplitude $[-n,0]$ (see \cite[\S 2.4]{tova}). We use the
same notation for $\dSch_{A_{i}}^{\leq n} \subset \dSch_{A_{i}}$. The base change
$\s$-functors
$$\dSch_{A_{i}} \longrightarrow \dSch_{A}$$
preserves cotangent complexes (see \cite[Lem. 1.4.1.16 (2)]{hagII}), and thus restrict to $\s$-functors
$$\dSch_{A_{i}}^{\leq n} \longrightarrow \dSch^{\leq n}_{A}.$$

\begin{prop}\label{p3}
The $\s$-functor
$$Colim\, \dSch_{A_{i}}^{\leq n} \longrightarrow \dSch_{A}^{\leq n}$$
is an equivalence of $\s$-categories. 
\end{prop}

\textit{Proof:} The proof again goes along the same general lines as for the proofs
of propositions \ref{p1} and \ref{p2}. We start by considering 
$\dAff_{A}^{\leq n} \subset \dSch_{A}^{\leq n}$ the full sub $\s$-category 
consisting of affine objects. We define in the same way $\dAff_{A_{i}}^{\leq n} \subset \dSch_{A_{i}}^{\leq n}$. 

\begin{lem}\label{l1}
The $\s$-functor
$$Colim\, \dAff_{A_{i}}^{\leq n} \longrightarrow \dAff_{A}^{\leq n}$$
is an equivalence. 
\end{lem}

\textit{Proof of the lemma:} We start by proving that the $\s$-functor is fully faithful. 
Let $B$ and $C$ be two simplicial $A_{i}$-algebras of finite presentation (and with
cotangent complexes relative to $A_{i}$ of amplitude $[-n,0]$). As $B$ is a finitely presented
$A_{i}$-algebra, we have
$$Map_{A}(B\otimes_{A_{i}}^{\mathbb{L}}A,C\otimes_{A_{i}}^{\mathbb{L}}A) \simeq
Map_{A_{i}}(B,Colim_{j\geq i}\, (C\otimes_{A_{i}}^{\mathbb{L}}A_{j}))\simeq
Colim_{j\geq i}\, Map_{A_{i}}(B,C\otimes_{A_{i}}^{\mathbb{L}}A_{j}),$$
where we have denoted by $Map_{A}$ the mapping spaces of the $\s$-category of
commutative simplicial $A$-algebras (and similarly for $A_{i}$). This proves fully faithfulness. 

Now, let $B$ a commutative simplicial $A$-algebra of finite presentation and with
cotangent complex $\mathbb{L}_{B/A}$ of amplitude in $[-n,0]$. We know that 
$B$ is equivalent to a retract of a finite cell commutative $A$-algebra $B'$. In particular
there is an index $i$ and a finite cell commutative $A_{i}$-algebra $B'_{i}$ such that 
$B' \simeq B_{i}'\otimes_{A_{i}}^{\mathbb{L}}A$. 

The object $B$ defines a projector up to homotopy on $B'$, that is a projector $p$
on $B'$ considered as an object in the homotopy category $Ho(A-CAlg)$ of 
commutative $A$-algebras. By chosing $i$ big enough we can 
moreover assume that this projector is induced by a projector $p_{i}$ on $B'_{i}$ in 
$Ho(A_{i}-sComm)$. By \cite[Sublemma 3]{to3} the category $Ho(A_{i}-sComm)$ is Karoubian closed, so 
$p_{i}$ splits as a composition in $Ho(A_{i}-sComm)$
$$p=vu : \xymatrix{ B'_{i} \ar[r]^-{u} & B_{i} \ar[r]^-{v} & B'_{i}}$$
for some commutative $A_{i}$-algebra $B_{i}$ and with $uv=id$. 

The $A_{i}$-algebra $B_{i}$ is of finite presentation (because it is a retract of
a finite cell commutative $A_{i}$-algebra) and we have
$B_{i} \otimes^{\mathbb{L}}_{A_{i}}A \simeq B$. 
It remains to show that $i$ can be chosen so that the cotangent complex
$\mathbb{L}_{B_{i}/A_{i}}$ has amplitude contained in $[-n,0]$. 
We let $Z_{i}:=Spec\, B_{i}$, $X_{i}:=Spec\, A_{i}$, $Z:=Spec\, B$ and $X=Spec\, A$. The locus in $Z_{i}$ in which 
$\mathbb{L}_{B_{i}/A_{i}}$ has amplitude contained in $[-n,0]$ is an
open  derived sub-scheme $U_{i} \subset Z_{i}$. Moreover, we have
$U_{i}\times_{X_{i}}X_{j} \simeq U_{j}$ for all $j\geq i$, because cotangent complexes
are stable by base changes. As 
$U_{i}\times_{X_{i}}X\simeq X$ there is an index $j$ with $U_{i}\times_{X_{i}}X_{j} \simeq X_{j}$
(we can cover $U_{i}$ by elementary opens $Spec\, B_{i}[f_{\alpha}^{-1}]$, 
and as $1$ is a linear combination of the $f_{\alpha}'s$ in $\pi_{0}(B)$ it must be so 
in some $\pi_{0}(B_{j})$). This finishes the proof of the lemma. \hfill $\Box$ \\
 
\begin{lem}\label{l2}
Let $X_{i}, Y_{i} \in \dSch_{A_{i}}^{\leq n}$, and denote by $X_{j}:=X_{i}\times_{Spec\, A_{i}}Spec\, A_{j}$ and
$Y_{j}:=Y_{i}\times_{Spec\, A_{i}}Spec\, A_{j}$ for $j\geq i$, and $X:=X_{i}\times_{Spec\, A_{i}}Spec\, A$, 
$Y:=Y_{i}\times_{Spec\, A_{i}}Spec\, A$.
Then, the natural morphism
$$Colim_{j\geq i}\, Map_{\dSch_{A_{j}}}(X_{j},Y_{j}) \longrightarrow Map_{\dSch_{A}}(X,Y)$$
is an equivalence. 
\end{lem}

\textit{Proof of the lemma:} Let us first assume that $X_{i}=Spec\, B_{i}$ is affine. 
We have $X=Spec\, B$ with $B=Colim\, B_{i}$. We thus have
$$Map_{\dSch_{A}}(X,Y) \simeq Map_{\dSch_{A_{i}}}(X,Y_{i}) \simeq
Colim_{j\geq i}\, Map_{\dSch_{A_{i}}}(X_{j},Y_{i}) \simeq 
Colim_{j\geq i}\, Map_{\dSch_{A_{j}}}(X_{j},Y_{j})$$
because $Y$ is locally of finite presentation over $Spec\, A$. To pass from the case where $X_{i}$ is affine
to the general case we use the same argument as for the proof of propositions \ref{p1} and \ref{p2}. 
The full sub $\s$-category of $\dSch_{A_{i}}^{\leq n}$ for which the lemma is true contains affine
and is stable by finite strong colimits. Therefore it contains all quasi-compact and quasi-separated derived
schemes. \hfill $\Box$ \\

To finish proposition \ref{p3} it remains to prove essential surjectivity. 

\begin{lem}\label{l3}
Let $f_{i} : X_{i} \longrightarrow Y_{i}$ be a morphism in $\dSch_{A_{i}}^{\leq n}$. Let 
$f : X \longrightarrow Y$ be the induced morphism by base change along $Spec\, A \longrightarrow Spec\, A_{i}$. 
If $f$ is a Zariski open immersion then there is $j\geq i$ for which 
$$f_{j}  : X_{j}=X_{i}\times_{Spec\, A_{i}}Spec\, A_{j} \longrightarrow 
Y_{j}=X_{i}\times_{Spec\, A_{i}}Spec\, A_{j}$$
is so.
\end{lem}

\textit{Proof of the lemma:} 
As $f$ is a Zariski open immersion it is
an \'etale monomorphism. Therefore, the cotangent complex $\mathbb{L}_{f}$ is zero. 
By compatibility of cotangent complexes by base changes, and by the proposition \ref{p2} we must have 
$\mathbb{L}_{f_{j}}$ for some $j\geq i$. Therefore, there is $j\geq i$ such that $f_{j}$ is \'etale. 
Moreover, as $f$ is a monomorphism the diagonal morphism $X \longrightarrow X\times_{Y}X$ is an equivalence. 
By the lemma \ref{l2} we must have a $j\geq i$ such that 
$X_{j} \longrightarrow X_{j}\times_{Y_{j}}X_{j}$ is an equivalence, or in other words that 
$f_{j}$ is a monomorphism of derived schemes. Therefore, there is a $j\geq i$ such that 
$f_{j}$ is an \'etale monomorphism and thus an open immersion. \hfill $\Box$ \\

We finally finish the proof of the proposition \ref{p3}. By lemma \ref{l1} we know that affine derived schemes
in $\dSch_{A}^{\leq n}$ belongs to the essential image. We first extend this to any quasi-compact quasi-affine derived scheme $X$ (i.e. any quasi-compact open of a derived affine scheme). 
Indeed, any such object is the image of 
a finite family of Zariski open affines
$$\{Spec\, B_{\alpha} \subset Spec\, B\}_{\alpha}.$$
By the lemmas \ref{l1} and \ref{l3}, this family is induced by a finite family of opens
$$\{Spec\, B_{i,\alpha} \subset Spec\, B_{i}\}_{\alpha},$$
in $\dSch_{A_{i}}^{\leq n}$ for some $i$. The image of this family defines 
a derived scheme $X_{i} \in \dSch_{A_{i}}^{\leq n}$ such that 
$X_{i}\times_{Spec\, A_{i}}Spec\, A \simeq X$. 

Finally, we proceed by induction on the number of  affines in an open covering. We assume that we have proven 
that all $X\in \dSch^{\leq n}_{A}$ covered by $k$ affine opens are in the essential image. If $X \in \dSch_{A}$
is covered by $(k+1)$ affine opens, we can form a push-out square in $\dSch_{A}^{\leq n}$
$$\xymatrix{
W \ar[r] \ar[d] & U \ar[d] \\
V \ar[r] & X,}$$
where all morphisms are Zariski open immersions, $V$ is affine and $U$ can be covered by $k$ affine opens. 
By what we have seen for quasi-affines, by induction and by lemmas \ref{l1} and \ref{l3}, this push-out square 
descent to a diagram of open immersions in $\dSch_{A_{i}}^{\leq n}$ for some $i$
$$\xymatrix{
W_{i} \ar[r] \ar[d] & U_{i} \\
V_{i}.  & }$$
Taking the push-out in $\dSch_{A_{i}}$ defines a derived scheme $X_{i} \in \dSch_{A_{i}}^{\leq n}$
such that $X_{i}\times_{Spec\, A_{i}}Spec\, A\simeq X$. \hfill $\Box$ \\

\section{Proof of the main theorem}

We are now ready to prove the main theorem of this work.

\begin{thm}\label{t1}
Let $f : X \longrightarrow Y$ be a local complete intersection and proper morphism 
schemes, then $f_{*} : L_{qcoh}(X) \longrightarrow L_{qcoh}(Y)$ preserves
perfect complexes. 
\end{thm}

\textit{Proof:} The statement is local on the Zariski topology of $Y$, so we can assume
that $Y=Spec\, A$ is affine (note that when doing so $X$ and $Y$ become automatically quasi-compact
and separated).
We write $A=Colim\, A_{i}$, a filtered colimit such that 
$A_{i}$ is noetherian for all $i$. As $f$ is lci, $\mathbb{L}_{f}$ is perfect of amplitude $[-1,0]$, 
so $X$ lies in $\dSch_{A}^{\leq 1}$. By proposition 
\ref{p3} there is an index $i$ and an object $X_{i} \in \dSch_{A_{i}}^{\leq 1}$ with
$X_{i} \times_{Spec\, A_{i}}Spec\, A\simeq X$. On the level of truncations we have a fiber product 
in the category of schemes
$$h^{0}(X_{i}) \times_{Spec\, A_{i}}Spec\, A\simeq h^{0}(X).$$
We can thus use \cite[Thm. 8.10.5]{egaIV} to show that if $i$ is taken big enough the scheme 
$h^{0}(X_{i})$ is proper over $Spec\, A_{i}$. Therefore, $X_{i} \longrightarrow Spec\, A_{i}$
is a proper lci morphism of derived schemes. 

Let $E \in L_{parf}(X)$. By the proposition \ref{p2} we can chose $i$ big enough so that $E$ descend to 
$E_{i} \in L_{parf}(X_{i})$. By the proposition \ref{p1} the theorem will be proven
if we can prove that $\hh(X_{i},E_{i})$ is a perfect $A_{i}$-dg-module. 

\begin{lem}\label{l4}
Let $f : X \longrightarrow S=Spec\, A$ be a proper and lci morphism of derived schemes
with $A$ a noetherian ring. Then 
$$f_{*} : L_{qcoh}(X) \longrightarrow L_{qcoh}(S)$$
preserves perfect complexes.
\end{lem}

\textit{Proof of the lemma:} We will first need to recall the following local structure
theorem for derived schemes whose cotangent complexes have amplitude in $[-1,0]$ (also 
called \emph{quasi-smooth} in the litterature).

\begin{sublem}\label{l5}
With the same notation as above 
the quasi-coherent sheaves $h^{i}(X)$ are coherent on $h^{0}(X)$ and only 
a finite number of them are non-zero.
\end{sublem}

\textit{Proof of the sublemma:} This is a local statement so we can assume that $X=Spec\, B$
with $B$ a retract of a finite cell commutative simplicial $A$-algebra with
$\mathbb{L}_{B/A}$ of amplitude in $[-1,0]$. As the statement 
we would like to prove is stable by retracts, we can even assume that $B$
has a finite cell decomposition:
$$\xymatrix{
B_{0}=A \ar[r] & B_{(1)} \ar[r] & B_{(2)} \ar[r] & \dots \ar[r] & B_{(k)}=B.}$$
For all $i$, there is a push-out of commutative simplicial $A$-algebras
$$\xymatrix{
B_{(i)} \ar[r] & B_{(i+1)} \\
\otimes^{\alpha_{i+1}} A[\partial \Delta^{i+1}] \ar[u] \ar[r] & \otimes^{\alpha_{i+1}} A[\Delta^{i+1}], \ar[u] }$$
where $A[K]$ denotes the free commutative simplicial $A$-algebra generated by a simplicial
set $K$ and $\alpha_{i}$ is the number of $i$-dimensional cells. Let us consider the morphism $p : B_{(1)} 
\longrightarrow B$. This morphism
induces an isomorphism on $\pi_{0}$ and an epimorphism on $\pi_{1}$. Its (homotopy) fiber 
is therefore connected. Using a Postnikov decomposition of the morphism $p$, obstruction
theory (see \cite[Lem. 2.2.1.1]{hagII}), and the fact that $\mathbb{L}_{B/A}$ as amplitude $[-1,0]$, we see
that the morphism $p$ as a section up to homotopy. In other words, 
the derived scheme $Spec\, B$ is a retract of $Spec\, B_{(1)}$. Moreover, by definition of cell 
algebras, $Spec\, B_{(1)}$ sits into a cartesian square
$$\xymatrix{
Spec\, B_{(1)} \ar[r] \ar[d] & \mathbb{A}^{\alpha_{0}}_{A} \ar[d]_-{u} \\
\{0\} \ar[r] & \mathbb{A}^{\alpha_{1}}_{A},}$$
where the morphism $u$ is determined by the attaching map 
$$\otimes^{\alpha_{1}} A[\partial \Delta^{1}] \longrightarrow B_{(0)}=A[X_{1},\dots,X_{\alpha_{0}}].$$
This shows that $Spec\, B_{(1)}$ satisfies the conclusion of the sublemma, and thus so does $Spec\, B$. 
This finishes the proof of \ref{l5}
\hfill $\Box$ \\

The sublemma and the fact that $E$ is perfect implies that 
the sheaves $H^{i}(E)$ are easily seen to be coherent and only a finite number of them are non-zero. Therefore, 
by Grothendieck's finiteness theorem and d\'evissage we have that $f_{*}(E)$ is a bounded coherent complex on $S$.
Moreover, for a closed point $s : Spec\, k(s) \longrightarrow S$, the proposition \ref{p1} implies that 
$s^{*}(f_{*}(E))\simeq \hh(X_{s},E_{s})$, where $X_{s}$ is the fiber of $f$ at $s$ and $E_{s}$ 
the pull-back of $E$ on $X_{s}$. Another application of the sublemma 
implies that 
$h^{i}(X_{s})$ is coherent and non-zero only for a finite number of indices
$i$. Therefore, perfect complexes on $X_{s}$ are also with coherent and bounded cohomology sheaves. 
As $X_{s}$ is proper the complex of $k(s)$-vector spaces 
$\hh(X_{s},E_{s})$ is cohomologically bounded with finite dimensional cohomology. The fibers of
the bounded coherent complex
$f_*(E)$ at every closed point is thus
cohomologically bounded. It is therefore of finite Tor dimension and is thus a perfect complex on $S$.
\hfill $\Box$ \\

This finishes the proof of theorem \ref{t1}. \hfill $\Box$ \\

\end{document}